\documentclass[12pt]{article}
\usepackage[utf8]{inputenc}
\usepackage[T1]{fontenc}
\usepackage{amsmath,amssymb,amsthm}
\usepackage{geometry}
\usepackage{hyperref}
\usepackage{indentfirst}
\usepackage{booktabs}
\usepackage{xcolor}
\usepackage{fvextra}
\DefineVerbatimEnvironment{code}{Verbatim}{
  breaklines=true,
  breakanywhere=true,
  fontsize=\footnotesize,
  xleftmargin=1em,
}
\DefineVerbatimEnvironment{term}{Verbatim}{
  breaklines=true,
  breakanywhere=true,
  fontsize=\footnotesize,
  xleftmargin=1em,
}

\newtheorem{theorem}{Theorem}
\newtheorem{lemma}{Lemma}
\newtheorem{proposition}{Proposition}
\newtheorem{corollary}{Corollary}
\theoremstyle{remark}
\newtheorem{remark}{Remark}
\theoremstyle{definition}
\newtheorem{definition}{Definition}

\title{Irreducibility of the Cuboid Polynomial $Q_{p,q}(t)$ via the Rank-One Chabauty--Coleman Method}

\author{Valery Asiryan\\[3pt]
\small \texttt{asiryanvalery@gmail.com}}
\date{\small November 4, 2026}

\begin{document}

\maketitle

% ===================== Abstract & Keywords =====================

\begin{abstract}
Let $Q_{p,q}(t)\in\mathbb{Z}[t]$ be Sharipov's even monic degree-$10$ second cuboid polynomial, with coprime positive integers $p\ne q$.
We prove that $Q_{p,q}(t)$ is irreducible over $\mathbb{Q}$, and hence over $\mathbb{Z}$.
More generally, the normalized associated quintic $P_s(x)$ and its composition $P_s(t^2)$ are irreducible over $\mathbb{Q}$ for every $s\in\mathbb{Q}\setminus\{0,\pm1\}$.

The proof reduces linear and quadratic factors to rational points and horizontal conjugate pairs on one fixed genus-$2$ curve.
A certified rank upper bound for its Jacobian, together with an exact divisor calculation, gives a rank-one annihilating differential at the prime $3$.
For the quadratic-factor case, the symmetric coordinates of a horizontal pair lie on an explicit double cover of an elliptic quartic.
All rational affine points on this model are $3$-adically integral and reduce to six smooth residue classes.
The pulled-back differential is a unit in each class, proving that six displayed points form the complete list.
Their coordinates exclude every admissible parameter.
A negative real root of the irreducible quintic then supplies the additional degree argument required for irreducibility of the degree-$10$ polynomial.
The rank computation and all exact arithmetic certificates are specified explicitly for reproduction.

\medskip
\noindent{\bf Keywords:}
irreducibility; perfect cuboid; cuboid polynomials; hyperelliptic curves; Jacobians; Chabauty--Coleman method; rational divisors; computer-assisted proof.

\smallskip
\noindent{\bf Mathematics Subject Classification:} 11D41, 11G30, 14H25, 12E05, 11Y16.
\end{abstract}

% ===================== Introduction =====================

\section{Introduction}\label{sec:intro}

The perfect cuboid problem asks for a rectangular box with positive integer edges whose three face diagonals and space diagonal are integers.
In Sharipov's polynomial formulation~\cite{SharipovCuboids,SharipovSecond}, one encounters explicit even polynomials depending on integer parameters.
The second cuboid polynomial $Q_{p,q}(t)$ has degree $10$ and is conjectured to be irreducible for coprime positive integers $p\ne q$.

\medskip
\noindent\textbf{Goal.}
We prove this irreducibility statement.
The intermediate result is stronger than required by the original parametrization: the associated normalized quintic $P_s(x)$ is irreducible for every $s\in\mathbb{Q}_{>0}\setminus\{1\}$, without assuming that $s$ is a rational square.

\medskip
\noindent\textbf{Method.}
A reducible quintic has an irreducible factor of degree at most two.
An explicit transformation sends a root of $P_s$ to a fixed genus-$2$ curve $C$, with second coordinate $(s-1)/(s+1)$.
A rational root gives a rational point of $C$.
An irreducible quadratic factor gives a conjugate pair on $C$ whose two points have the same second coordinate.
Thus both cases are controlled by the same Jacobian.

We use the rank-one form of the Chabauty--Coleman method~\cite{Chabauty,Coleman}.
The treatment of conjugate pairs is related to symmetric-power methods~\cite{Siksek}, but only a one-dimensional locus of horizontal pairs is needed here.
We construct its model directly and prove the required local uniqueness rather than invoke a rational-point theorem for the entire symmetric square.
A rank bound for the genus-$2$ Jacobian is obtained computationally, in the setting of explicit descent~\cite{Stoll}, using \textsc{Magma}~\cite{Magma}.
All subsequent steps use a small exact divisor calculation, power series over $\mathbb{Z}_3$, and enumeration over $\mathbb{F}_3$.
In particular, the proof does not assume completeness of a bounded search for rational points.

The article is organized as follows.
Section~\ref{sec:poly} states the normalization and the main theorem.
Sections~\ref{sec:curve}--\ref{sec:annihilator} establish the rank-one differential.
Section~\ref{sec:linear} excludes linear factors.
Sections~\ref{sec:quadratic}--\ref{sec:six} determine the horizontal-pair obstruction and exclude quadratic factors.
Section~\ref{sec:conclusion} completes the degree-$10$ argument.
The appendix records the computational inputs and reproducible verification code.

% ===================== Polynomial and normalization =====================

\section{The second cuboid polynomial and normalization}\label{sec:poly}

Let $p,q\in\mathbb{Z}_{>0}$ be coprime and $p\ne q$.
Sharipov's second cuboid polynomial is
\begin{align}
Q_{p,q}(t)
={}&t^{10}+(2q^2+p^2)(3q^2-2p^2)t^8\nonumber\\
&+(q^8+10p^2q^6+4p^4q^4-14p^6q^2+p^8)t^6\nonumber\\
&-p^2q^2(q^8-14p^2q^6+4p^4q^4+10p^6q^2+p^8)t^4\nonumber\\
&-p^6q^6(q^2+2p^2)(-2q^2+3p^2)t^2-p^{10}q^{10}.
\label{eq:Qpq}
\end{align}
It is weighted-homogeneous of weight $20$ for
$\deg p=\deg q=1$ and $\deg t=2$.

\begin{definition}[Normalized second cuboid quintic]\label{def:Ps}
For $s\in\mathbb{Q}$, define
\begin{equation}\label{eq:Ps}
\begin{aligned}
P_s(x)={}&x^5+(2+s)(3-2s)x^4
 +(1+10s+4s^2-14s^3+s^4)x^3\\
&-s(1-14s+4s^2+10s^3+s^4)x^2
 -s^3(1+2s)(-2+3s)x-s^5.
\end{aligned}
\end{equation}
\end{definition}

\begin{lemma}[Normalization]\label{lem:normalization}
With $s=(p/q)^2$, one has
\begin{equation}\label{eq:normalization}
Q_{p,q}(t)=q^{20}P_s\!\left(\frac{t^2}{q^4}\right).
\end{equation}
\end{lemma}
\begin{proof}
Set $r=p/q$ and substitute $p=rq$, $t=q^2u$ into \eqref{eq:Qpq}.
Factoring out $q^{20}$ and setting $s=r^2$ gives $P_s(u^2)$.
\end{proof}

\begin{theorem}\label{thm:main}
For every $s\in\mathbb{Q}_{>0}\setminus\{1\}$, the quintic $P_s(x)$ is irreducible over $\mathbb{Q}$.
Consequently, for coprime positive integers $p\ne q$, the polynomial $Q_{p,q}(t)$ is irreducible over $\mathbb{Q}$ and over $\mathbb{Z}$.
\end{theorem}

The last implication requires a separate argument: irreducibility of $P_s(x)$ alone does not automatically give irreducibility of $P_s(t^2)$.
That argument will be supplied in Section~\ref{sec:conclusion}.

% ===================== Fixed curve and rank =====================

\section{A fixed genus-two curve and a rank bound}\label{sec:curve}

Put
\begin{align}
f(z)&=-16z^5+48z^4-48z^3+16z^2+z,\label{eq:f}\\
h(z)&=1-f(z)
 =(1-z)(1-2z)^2(1+4z-4z^2).\label{eq:h}
\end{align}
Let $C$ denote the smooth projective curve associated with
\begin{equation}\label{eq:C}
C:\qquad y^2=f(z),
\end{equation}
and write $J=\operatorname{Jac}(C)$.
The discriminant of $f$ is
\begin{equation}\label{eq:discriminant}
\operatorname{disc}(f)=-2^{16}\cdot43.
\end{equation}
Hence $C$ has genus $2$ and good reduction at $3$.
It has a unique rational point $\infty$ at infinity.
We also use the rational Weierstrass point $T=(0,0)$.
A basis of holomorphic differentials is
\begin{equation}\label{eq:omegas}
\omega_0=\frac{dz}{2y},\qquad
\omega_1=\frac{z\,dz}{2y}.
\end{equation}
Both are integral differentials at the prime $3$.

\begin{lemma}[Computational rank bound]\label{lem:rank}
The Mordell--Weil rank of $J(\mathbb{Q})$ is at most one.
\end{lemma}
\begin{proof}[Computer-assisted proof]
Consider the genus-$2$ model
\begin{equation}\label{eq:C0}
C_0:\qquad w^2=g(t):=t^5+21t^4+26t^3+10t^2+5t+1.
\end{equation}
The substitution
\begin{equation}\label{eq:isomorphism}
z=\frac{1}{t+1},\qquad y=\frac{w}{(t+1)^3}
\end{equation}
satisfies
\[
(t+1)^6 f\!\left(\frac{1}{t+1}\right)=g(t).
\]
Its inverse is $t=z^{-1}-1$, $w=y/z^3$ on the corresponding affine open sets.
Thus it extends to an isomorphism of smooth projective curves over $\mathbb{Q}$, and hence of their Jacobians.

Let $J_0=\operatorname{Jac}(C_0)$.
A certified computation in \textsc{Magma} gives
\[
\operatorname{rank} J_0(\mathbb{Q})\le1.
\]
The class-group calculations use full Minkowski bounds, so this upper bound is unconditional and does not rely on the generalized Riemann hypothesis.
The complete computation and its output are recorded in Appendix~\ref{app:scripts}.
Since \eqref{eq:isomorphism} induces an isomorphism $J\simeq J_0$ over $\mathbb{Q}$, their Mordell--Weil ranks are equal, proving the assertion.
\end{proof}

\begin{remark}\label{rem:rankrole}
Lemma~\ref{lem:rank} is the computational input to the proof.
We do not need a complete Mordell--Weil basis or a saturation calculation.
An explicit divisor below will have nonzero $3$-adic logarithm, proving that the rank is exactly one and determining its annihilating differential to the required precision.
\end{remark}

% ===================== Exact logarithm =====================

\section{An exact rank-one differential at the prime three}\label{sec:annihilator}

For a holomorphic differential $\omega$, let $\log_\omega$ denote its Abelian logarithm on $J(\mathbb{Q}_3)$; see~\cite{ColemanIntegrals} for the integration framework.
On a residue disk this logarithm is obtained by integrating the corresponding local power series.
The logarithm kills torsion and is additive on divisor classes.
Its compatibility with pullback will also be used for the horizontal-pair model.

\begin{lemma}[A local uniqueness criterion]\label{lem:local}
Let a residue disk have a parameter $q\in3\mathbb{Z}_3$, centered at $q=0$.
Suppose a differential on it is $b(q)\,dq$, where
$b(q)\in\mathbb{Z}_3[[q]]$ and $b(0)\in\mathbb{Z}_3^\times$.
Then the primitive vanishing at $q=0$ has no other zero in this disk.
\end{lemma}
\begin{proof}
Write $b(q)=\sum_{j\ge0}b_jq^j$.
For $q\ne0$, its primitive divided by $q$ is
\[
b_0+\sum_{j\ge1}\frac{b_jq^j}{j+1}.
\]
For each $j\ge1$, one has
$v_3(q^j/(j+1))\ge j-v_3(j+1)\ge1$.
Thus this expression is a unit, proving the claim.
\end{proof}

We represent divisor classes by Mumford coordinates $(U,V)$: here $U$ is monic, $\deg V<\deg U\le2$, and $U\mid f-V^2$.
The corresponding class is the effective divisor supported on $U(z)=0$, $y=V(z)$, minus $(\deg U)\infty$.
Such coordinates allow exact addition by Cantor's algorithm~\cite{Cantor}.

\begin{lemma}[Annihilating differential]\label{lem:annihilator}
There exists a differential annihilating $J(\mathbb{Q})$ of the form
\begin{equation}\label{eq:annihilator}
\omega=a\omega_0+\omega_1,\qquad
 a\in\mathbb{Z}_3,\qquad a\equiv3\pmod9.
\end{equation}
\end{lemma}
\begin{proof}
Let
\[
A=[(1/2,1)-\infty]\in J(\mathbb{Q}).
\]
Exact divisor addition gives the following small certificate:
\begin{equation}\label{eq:11A}
11A=(U,V),\quad
\begin{aligned}
U(z)&=z^2-\frac{69}{196}z-\frac9{392},\\
V(z)&=-\frac{4085}{2058}z-\frac{233}{1372}.
\end{aligned}
\end{equation}
Appendix~\ref{app:scripts} specifies checks for both the divisor-class equality and the divisibility $U\mid f-V^2$.
The latter condition alone would not certify that the class is $11A$.

Let $(z_i,y_i)$, $i=1,2$, be the two geometric points represented by \eqref{eq:11A}.
Since $y_i=V(z_i)$, their trace and norm are obtained directly from $U$:
\begin{equation}\label{eq:trace-norm}
\tau:=y_1+y_2=-\frac{139623}{134456},\qquad
\nu:=y_1y_2=\frac{7671}{134456}.
\end{equation}
Write $\tau=3r$ and $\nu=3n$.
Then
\begin{equation}\label{eq:rn}
r,n\in\mathbb{Z}_3^\times,\qquad r\equiv n\equiv2\pmod3.
\end{equation}
The nonleading coefficients of $U$ have valuations $1$ and $2$.
Its Newton polygon therefore gives $v_3(z_i)=1$.
From $f(z)=z+O(z^2)$ it follows that $v_3(y_i)=1/2$.
Both geometric points lie in the disk of $T$ over a finite extension of $\mathbb{Q}_3$.

Near $T$, the coordinate $y$ is a uniformizer and the equation $y^2=f(z)$ has inverse series
\begin{equation}\label{eq:inverse}
z=y^2-16y^4+560y^6+\cdots\in\mathbb{Z}_3[[y^2]].
\end{equation}
The integrality follows by formal inversion, because $f'(0)=1$.
Since $\operatorname{div}(z)=2T-2\infty$, the class $[T-\infty]$ is $2$-torsion and $[2T-2\infty]=0$.
Consequently, with $S_m=y_1^m+y_2^m$, the logarithms of $11A$ are sums of integrals from $T$:
\begin{align}
\ell_0:=\log_{\omega_0}(11A)
 &=\tau-\frac{32}{3}S_3+\cdots,\label{eq:ell0}\\
\ell_1:=\log_{\omega_1}(11A)
 &=\frac{1}{3}S_3-\frac{48}{5}S_5+\cdots.\label{eq:ell1}
\end{align}
These formulas result by differentiating \eqref{eq:inverse}, substituting in \eqref{eq:omegas}, and integrating term by term.
Both local integrals are computed in a finite extension $K/\mathbb{Q}_3$ containing the two support points, with $v_3(3)=1$.
Since the unintegrated coefficients are integral and $v_3(y_i)=1/2$, a term of degree $2j-1$ in either primitive has valuation at least $j-\tfrac12-v_3(2j-1)$, which tends to infinity.
Thus the individual series converge in $K$ and may be added term by term to obtain the expressions in $S_{2j-1}$.

We make the precision bounds explicit.
The power sums satisfy $S_0=2$, $S_1=\tau$, and
\[
S_m=\tau S_{m-1}-\nu S_{m-2}.
\]
For $m=2j-1$, each monomial in the expression for $S_m$ in terms of $\tau,\nu$ contains at least $j$ factors, so
\begin{equation}\label{eq:powerbound}
v_3(S_{2j-1})\ge j.
\end{equation}
Moreover,
\begin{equation}\label{eq:S3}
S_3=\tau^3-3\nu\tau=27r(r^2-n).
\end{equation}
All unintegrated power-series coefficients are integral; integration introduces only the denominator $2j-1$ in the term involving $S_{2j-1}$.
In \eqref{eq:ell0}, the $S_3$ term therefore has valuation at least $2$, and every subsequent term also has valuation at least $2$.
In \eqref{eq:ell1}, all terms after $S_3$ have valuation at least $3$, since
\[
j-v_3(2j-1)\ge3\qquad(j\ge3).
\]
The valuations tend to infinity, so these estimates also justify passage to the infinite sums.
It follows that
\begin{equation}\label{eq:logcongruences}
\ell_0\equiv3r\pmod9,\qquad
\ell_1\equiv9r(r^2-n)\pmod{27}.
\end{equation}
In particular, $\ell_0\ne0$, so $A$ has infinite order.
By Lemma~\ref{lem:rank}, the group $J(\mathbb{Q})$ has rank exactly one.
Its logarithmic image spans the same $\mathbb{Q}_3$-line as $\log(11A)$.

Set $a=-\ell_1/\ell_0$.
This gives an annihilating differential for every rational divisor class.
The congruences above imply $a\in3\mathbb{Z}_3$ and
\[
\frac{a}{3}\equiv-(r^2-n)\equiv1\pmod3.
\]
Thus $a\equiv3\pmod9$, as required.
The argument only uses a nonzero vector in the rank-one logarithmic image; no assertion that $A$ generates the full Mordell--Weil group is needed.
\end{proof}

% ===================== Linear factors =====================

\section{Rational points and exclusion of linear factors}\label{sec:linear}

\begin{proposition}[Rational points on $C$]\label{prop:Cpoints}
The full rational-point set is
\begin{equation}\label{eq:Cpoints}
C(\mathbb{Q})=\{\infty,(0,0),(1,\pm1),(1/2,\pm1)\}.
\end{equation}
\end{proposition}
\begin{proof}
Enumeration over $\mathbb{F}_3$ gives
\begin{equation}\label{eq:Cmod3}
C(\mathbb{F}_3)=\{\infty,(0,0),(1,\pm1),(2,\pm1)\}.
\end{equation}
Every rational point reduces to one of these six points by properness of the smooth model.
At each of the four finite points with nonzero first coordinate, the reduction of $\omega$ is the nonzero differential $z\,dz/(2y)$.
It is also nonzero at infinity.
Indeed, the local parameter $z^2/y$ at infinity shows that $z\,dz/(2y)$ has unit leading coefficient, whereas $dz/(2y)$ vanishes there.
Each of these five disks contains exactly one of the displayed rational points.
The Abelian integral of $\omega$ between rational points vanishes, so Lemma~\ref{lem:local} gives uniqueness in each disk.

It remains to examine the disk of $T=(0,0)$, where the reduced differential vanishes.
Use $y\in3\mathbb{Z}_3$ as parameter and put
\[
F_T(y)=\int_T^{(z(y),y)}\omega.
\]
From \eqref{eq:inverse} one obtains, for $y\ne0$,
\begin{equation}\label{eq:Tprimitive}
\frac{F_T(y)}{y}
=a+\frac{1-32a}{3}y^2+R(y),\qquad v_3(R(y))\ge4.
\end{equation}
For the remainder bound, the terms remaining after division by $y$ have degree $2k$, $k\ge2$, and their coefficients belong to $(2k+1)^{-1}\mathbb{Z}_3$.
For $y\in3\mathbb{Z}_3$ use $2k-v_3(2k+1)\ge4$ to obtain the stated valuation bound.

If $v_3(y)=1$, write $y=3c$ with $c\in\mathbb{Z}_3^\times$.
Using $a/3\equiv1\pmod3$ gives
\[
\frac{F_T(y)}{3y}\equiv1+c^2\equiv2\pmod3.
\]
If $v_3(y)\ge2$, the term $a$ on the right of \eqref{eq:Tprimitive} has valuation one, whereas every other term has larger valuation.
Hence $F_T(y)\ne0$ whenever $y\ne0$.
Since the integral must vanish at any rational point in this disk, $T$ is its only rational point.
This completes the list.
\end{proof}

We now connect this fixed curve to the parameter-dependent quintic.
For $s\ne0$, define the reciprocal polynomial
\begin{equation}\label{eq:R}
R_s(x)=-\frac{x^5}{s^5}P_s\!\left(\frac{s^2}{x}\right)\in\mathbb{Q}[x].
\end{equation}
Its leading coefficient is one.
Put
\begin{equation}\label{eq:rootmap}
\Delta_s(x)=(x+1)(x+s^2),\qquad
Z_s(x)=\frac{x(s-1)^2}{\Delta_s(x)},\qquad
v_s=\frac{s-1}{s+1}.
\end{equation}
A direct polynomial identity is
\begin{equation}\label{eq:rootidentity}
(s+1)^2\Delta_s(x)^5\bigl(f(Z_s(x))-v_s^2\bigr)
=-(s-1)^2P_s(x)R_s(x).
\end{equation}
Both sides become polynomials in $s,x$ after clearing the displayed substitution; Appendix~\ref{app:scripts} includes its exact verification.

The transformation can also be motivated by inversion invariants.
With $\rho=x/s$, set $U=s+1/s$ and $V=\rho+1/\rho$.
The quotient parameter
\[
t=\frac{V+2}{U-2}=\frac{(x+s)^2}{x(s-1)^2}
\]
gives $z=1/(t+1)=Z_s(x)$, in agreement with \eqref{eq:isomorphism}.
For the proof itself, the identity \eqref{eq:rootidentity} is sufficient.

\begin{lemma}[Roots give points of $C$]\label{lem:rootmap}
For $s\in\mathbb{Q}_{>0}\setminus\{1\}$, every root $x$ of $P_s$ gives a point
$(Z_s(x),v_s)$ on $C$, defined over $\mathbb{Q}(x)$.
Its first coordinate is finite and nonzero.
\end{lemma}
\begin{proof}
The relevant exceptional values are
\begin{align}
P_s(0)&=-s^5,\nonumber\\
P_s(-1)&=4(s-1)^4(s+1),\label{eq:specialvalues}\\
P_s(-s^2)&=-4s^5(s-1)^4(s+1).\nonumber
\end{align}
They are nonzero under the stated assumptions.
Thus $Z_s(x)$ is defined and nonzero at every root.
The identity \eqref{eq:rootidentity} now gives $f(Z_s(x))=v_s^2$.
\end{proof}

\begin{corollary}[No linear factor]\label{cor:linear}
For $s\in\mathbb{Q}_{>0}\setminus\{1\}$, the polynomial $P_s$ has no rational root.
\end{corollary}
\begin{proof}
A rational root would give a rational point of $C$ whose second coordinate is $v_s$.
But $-1<v_s<1$ and $v_s\ne0$, whereas Proposition~\ref{prop:Cpoints} permits only second coordinates $0,1,-1$ at finite rational points.
\end{proof}

% ===================== Horizontal pair model =====================

\section{Quadratic factors and horizontal conjugate pairs}\label{sec:quadratic}

We first check that distinct roots of $P_s$ cannot coalesce under the map $Z_s$.
This matters because dividing a difference of curve equations by the difference of their first coordinates will be essential.

\begin{lemma}[No collision of root images]\label{lem:collision}
For $s\in\mathbb{Q}_{>0}\setminus\{1\}$, two distinct roots $x_1,x_2$ of $P_s$ have distinct images under $Z_s$.
\end{lemma}
\begin{proof}
Coefficient comparison in \eqref{eq:Ps} and \eqref{eq:R} gives
\begin{equation}\label{eq:reciprocal}
P_s(x)-R_s(x)=8(1-s^2)x(x+s)^3,
\qquad P_s(-s)=-2s^3(s-1)^4.
\end{equation}
Thus $P_s$ cannot have both $x$ and $s^2/x$ as roots.
Indeed, such an $x$ would be a common root of $P_s,R_s$, forcing $x=0$ or $x=-s$, and neither is a root.
On the other hand,
\begin{equation}\label{eq:zdifference}
Z_s(x_1)-Z_s(x_2)=
\frac{(s-1)^2(x_1-x_2)(s^2-x_1x_2)}
{\Delta_s(x_1)\Delta_s(x_2)}.
\end{equation}
All denominators are nonzero by \eqref{eq:specialvalues}.
Equality of the images would therefore force $x_1x_2=s^2$, contrary to the preceding observation.
\end{proof}

Suppose $P_s$ has a quadratic divisor over $\mathbb{Q}$.
By Corollary~\ref{cor:linear} this divisor is irreducible, so its roots $x_1,x_2$ are distinct quadratic conjugates.
By Lemmas~\ref{lem:rootmap} and~\ref{lem:collision}, their images
\[
(z_1,v_s),\ (z_2,v_s)\in C(\overline{\mathbb{Q}})
\]
are distinct conjugate points.
Their common second coordinate makes this a horizontal conjugate pair.
Write
\begin{equation}\label{eq:symmetric}
e=z_1+z_2,\qquad m=z_1z_2,\qquad
u=e-1,\qquad Y=4m-6u^2.
\end{equation}
All four quantities are rational.
Since $h(z_1)=h(z_2)$, division by $z_1-z_2$ gives
\begin{equation}\label{eq:em}
16m^2-48(e-1)^2m+16e(e-1)^3-1=0.
\end{equation}
After completing the square in $m$, this becomes the quartic
\begin{equation}\label{eq:E}
E:\qquad Y^2=20u^4-16u^3+1.
\end{equation}
The smooth projective completion of this nonsingular quartic, with the rational point $(0,1)$ chosen as origin, is an elliptic curve.

To recover the common value of $h$, reduce it modulo $z^2-ez+m$.
Using \eqref{eq:em}, the remainder is constant and equals
\begin{align}
W&=1-16mu^3+16m^2(2u-1)\nonumber\\
 &=2u\bigl(10Yu^2-6Yu+44u^4-44u^3+8u^2+1\bigr),
\label{eq:W}
\end{align}
where $m=(Y+6u^2)/4$ in the second line.
Since $f=1-h$, the horizontal pair gives a rational point of the affine model
\begin{equation}\label{eq:X}
X:\qquad
Y^2=20u^4-16u^3+1,\qquad v^2=1-W(u,Y),
\end{equation}
with $v=v_s$.
Thus the quadratic-factor problem has become a rational-point problem on a fixed curve.

\begin{lemma}[Divisor map and pulled-back differentials]\label{lem:pullback}
For a point of $X$, put
\[
e=u+1,\qquad m=\frac{Y+6u^2}{4},\qquad B(z)=z^2-ez+m.
\]
The two roots of $B$ define an effective degree-two divisor
\[
D=(z_1,v)+(z_2,v)
\]
on $C$, with multiplicities when the roots coincide.
The class $\Phi=[D-2\infty]$ is rational at rational points of $X$.
On the locus $v\ne0$, the pullbacks of \eqref{eq:omegas} are
\begin{equation}\label{eq:pullback}
\eta_0=\frac{du}{2v},\qquad
\eta_1=\frac{(1-2u)\,du-\tfrac14dY}{2v}.
\end{equation}
\end{lemma}
\begin{proof}
The reduction calculation above gives $B\mid h-W$.
By \eqref{eq:X}, it follows that $B\mid f-v^2$, so the indicated horizontal divisor lies on $C$.
Its definition is Galois invariant when $u,Y,v$ are rational.
The map to the Jacobian is given by its degree-zero divisor class.

For distinct roots of $B$, the pullback of a differential to this class is the sum of its values on the two points.
Since $z_1+z_2=e$ and $z_1z_2=m$,
\[
dz_1+dz_2=de=du,
\qquad z_1dz_1+z_2dz_2=e\,de-dm.
\]
Now $dm=\tfrac14dY+3u\,du$, which gives \eqref{eq:pullback}.
The formulas are symmetric in the roots and extend across coincident roots wherever $v\ne0$.
\end{proof}

The differential
\begin{equation}\label{eq:eta}
\eta=a\eta_0+\eta_1
\end{equation}
is the pullback of the annihilator in Lemma~\ref{lem:annihilator}.
For rational points $R,R_0$ of $X$ in a common smooth residue disk on which these expressions are integral, its local integral satisfies
\begin{equation}\label{eq:integralzero}
\int_{R_0}^{R}\eta
=\log_\omega\bigl(\Phi(R)-\Phi(R_0)\bigr)=0.
\end{equation}
In the disks used below, $Y$ and $v$ are units and all denominators are units at $3$.
The divisor construction and the displayed differentials are therefore defined integrally on each disk.

% ===================== Six disks =====================

\section{Six residue disks determine the obstruction points}\label{sec:six}

\begin{proposition}[Complete affine rational-point list]\label{prop:Xpoints}
The rational affine points of \eqref{eq:X} are exactly
\begin{equation}\label{eq:Xpoints}
(u,Y,v)=(0,1,\pm1),\qquad(0,-1,\pm1),\qquad(1/2,1/2,\pm1).
\end{equation}
\end{proposition}
\begin{proof}
First, every rational affine point is integral at $3$.
If $v_3(u)<0$, divide \eqref{eq:E} by $u^4$:
\[
\left(\frac{Y}{u^2}\right)^2=20-\frac{16}{u}+\frac1{u^4}.
\]
The right side is a unit with reduction $2$ in $\mathbb{F}_3$, which cannot be the reduction of a square in $\mathbb{Q}_3$.
Hence $u\in\mathbb{Z}_3$.
The two equations of $X$ then give $Y,v\in\mathbb{Z}_3$.
Consequently every rational affine point has a reduction in this affine model.

Enumerating the equations over $\mathbb{F}_3$ yields precisely
\begin{equation}\label{eq:Xmod3}
(u,Y,v)=(0,1,\pm1),\qquad(0,2,\pm1),\qquad(2,2,\pm1).
\end{equation}
At all six points both $Y$ and $v$ are nonzero.
The derivatives of the two equations with respect to $Y,v$ form an invertible triangular matrix, with diagonal entries $2Y$ and $2v$.
Thus all six points are smooth, and $u-u_0$ is a local parameter in the disk of any chosen lift $(u_0,Y_0,v_0)$.
Each disk contains one of the six rational points in \eqref{eq:Xpoints}; direct substitution gives $W=0$ at those points.

By Lemma~\ref{lem:annihilator}, $a\equiv0\pmod3$, so $\eta$ reduces to $\eta_1$.
Differentiating \eqref{eq:E} gives
\[
dY=\frac{80u^3-48u^2}{2Y}\,du.
\]
Therefore
\begin{equation}\label{eq:etacoefficient}
\eta_1=
\frac{1-2u-(80u^3-48u^2)/(8Y)}{2v}\,du.
\end{equation}
The numerator reduces to $1$ at either point with $u=0$, and to $2$ at $(u,Y)=(2,2)$.
Since $v$ is a unit, the coefficient of $du$ is a unit in every disk in \eqref{eq:Xmod3}.

Fix the known rational point $R_0$ in one of these disks.
Any other rational point $R$ in the disk satisfies \eqref{eq:integralzero}.
Lemma~\ref{lem:local}, applied with parameter $u-u(R_0)$, forces $R=R_0$.
Every rational affine point belongs to one of the six disks, so the displayed list is complete.
\end{proof}

\begin{corollary}[No quadratic factor]\label{cor:quadratic}
For $s\in\mathbb{Q}_{>0}\setminus\{1\}$, the polynomial $P_s$ has no quadratic divisor over $\mathbb{Q}$.
\end{corollary}
\begin{proof}
The construction of Section~\ref{sec:quadratic} sends any such divisor to a rational affine point of $X$ with $v=(s-1)/(s+1)$.
For positive $s$ this satisfies $-1<v<1$.
Proposition~\ref{prop:Xpoints} gives $v=\pm1$ at every rational affine point of $X$, a contradiction.
\end{proof}

\begin{remark}
The integrality argument is essential: it ensures that the six affine residue disks cover every rational point required in the factorization argument.
No bound on the numerator or denominator of $u$ is used.
The uniqueness in each disk then proves completeness, rather than merely confirming a list found by search.
\end{remark}

% ===================== Final degree argument =====================

\section{Completion of the irreducibility proof}\label{sec:conclusion}

\begin{proof}[Proof of Theorem~\ref{thm:main}]
A reducible polynomial of degree five over a field has an irreducible factor of degree at most two.
Indeed, if every irreducible factor had degree at least three, two nonconstant factors would already have total degree at least six.
Corollaries~\ref{cor:linear} and~\ref{cor:quadratic} therefore prove irreducibility of $P_s(x)$ over $\mathbb{Q}$.

We next establish irreducibility after composition with a square.
By \eqref{eq:specialvalues},
\[
P_s(-1)=4(s-1)^4(s+1)>0,\qquad P_s(0)=-s^5<0.
\]
Thus $P_s$ has a real root $\alpha\in(-1,0)$.
Its irreducibility gives $[\mathbb{Q}(\alpha):\mathbb{Q}]=5$.
View $\mathbb{Q}(\alpha)$ as the subfield of $\mathbb{R}$ defined by this negative real root.
A square in this real field is nonnegative, so $\alpha$ is not a square in $\mathbb{Q}(\alpha)$.
The polynomial $u^2-\alpha$ is consequently irreducible over $\mathbb{Q}(\alpha)$.

If $\beta^2=\alpha$, then $\mathbb{Q}(\alpha)\subset\mathbb{Q}(\beta)$ and
\[
[\mathbb{Q}(\beta):\mathbb{Q}]
=[\mathbb{Q}(\beta):\mathbb{Q}(\alpha)]\,[\mathbb{Q}(\alpha):\mathbb{Q}]
=2\cdot5=10.
\]
Since $\beta$ is a root of the monic degree-$10$ polynomial $P_s(u^2)$, that polynomial is irreducible over $\mathbb{Q}$.
For $s=(p/q)^2$, the invertible rational change of variable $u=t/q^2$ and the nonzero scalar factor in \eqref{eq:normalization} preserve irreducibility.
Hence $Q_{p,q}(t)$ is irreducible over $\mathbb{Q}$.
Finally, it is monic with integer coefficients, so Gauss's lemma gives irreducibility over $\mathbb{Z}$.
\end{proof}

\begin{corollary}[Classification of rational parameters]\label{cor:rational-parameters}
For $s\in\mathbb{Q}$, both $P_s(x)$ and $P_s(t^2)$ are irreducible over $\mathbb{Q}$ if and only if $s\notin\{0,\pm1\}$.
\end{corollary}
\begin{proof}
For $s\notin\{0,\pm1\}$, the proofs of Lemmas~\ref{lem:rootmap} and~\ref{lem:collision} apply unchanged: all exceptional values and denominators used there remain nonzero.
Moreover, $v_s=(s-1)/(s+1)\notin\{0,\pm1\}$.
Propositions~\ref{prop:Cpoints} and~\ref{prop:Xpoints} therefore exclude linear and quadratic factors of $P_s$, proving its irreducibility.

For $s>0$, irreducibility of $P_s(t^2)$ follows from the degree argument above.
For $s<0$, we have $P_s(0)>0$ and $P_s(x)\to-\infty$ as $x\to-\infty$, so a negative real root exists and the same argument applies.
Finally,
\[
\begin{aligned}
P_0(x)&=x^3(x^2+6x+1),\\
P_1(x)&=(x-1)(x+1)^4,\\
P_{-1}(x)&=(x+1)^5,
\end{aligned}
\]
which proves reducibility in all exceptional cases, also after substituting $x=t^2$.
\end{proof}

% ===================== Appendix =====================

\appendix
\section{\textsc{Magma} scripts and transcripts}\label{app:scripts}

The only rank computation used in the proof is the upper bound in Lemma~\ref{lem:rank}.
Except for the auxiliary genus-two rational-point verification, the remaining computational checks use exact rational polynomial arithmetic or finite enumeration over $\mathbb{F}_3$.
The code below is arranged to keep these two roles separate.
Its rank target is $1$; the remaining assertions check the isomorphism of curve models, the root-map and reciprocal identities, the divisor certificate, its trace and norm, and the residue-class calculations.
The analytic uniqueness arguments are proved in the main text and do not depend on a rational-point search.

\subsection*{Script: Complete exact verification}
The script was executed with \textsc{Magma} V2.29-10.
The two class-group bound maps are both set to the Minkowski bound, so the rank calculation is unconditional.
The additional \texttt{Chabauty} call independently reproduces the genus-two rational-point list; the proofs in the main text require only the rank upper bound and the exact local certificate.
The identities and finite-field checks are assertions, so an unsuccessful check must be resolved before treating the run as a certificate.

In the final part of the script, the variables \texttt{e,f,w} denote the symmetric coordinates $(e,m,v)$.
They are related to the coordinates of the main text by
\[
u=e-1,\qquad Y=4m-6(e-1)^2,\qquad v=w;
\]
the auxiliary variable \texttt{V} in that part equals $Y/2$.
This change of coordinates is invertible over $\mathbb{Z}_3$.
Consequently, the printed rational lifts and residue classes give exactly the lists in \eqref{eq:Xpoints} and \eqref{eq:Xmod3}, respectively, in symmetric coordinates.
The nonvanishing checks use these same coordinates and verify the unit condition for the pulled-back differential.

\newpage
\noindent\textbf{Code.}
\begin{code}
Q := Rationals();
Rt<t> := PolynomialRing(Q);
F0 := t^5+21*t^4+26*t^3+10*t^2+5*t+1;
C0 := HyperellipticCurve(F0);
J0 := Jacobian(C0);
proof_bound := map<PowerStructure(RngOrd) -> Integers() |
    O :-> MinkowskiBound(O)>;
SetClassGroupBoundMaps(proof_bound, proof_bound);
upper := RankBound(J0);
print "Unconditional rank upper bound:", upper;
assert upper eq 1;
// Independent exact verification of the genus-two rational-point list.
B0 := J0![C0![0,1,1], C0![1,0,0]];
assert Order(B0) eq 0;
points0 := Chabauty(B0);
expected0 := {C0![1,0,0], C0![-1,0,1], C0![0,1,1],
    C0![0,-1,1], C0![1,8,1], C0![1,-8,1]};
assert #points0 eq 6;
assert &and[p in points0 : p in expected0];
print "Complete genus-two rational-point list:", points0;
Rz<z> := PolynomialRing(Q);
F := -16*z^5+48*z^4-48*z^3+16*z^2+z;
assert Discriminant(F) eq -2^16*43;
Kz := FieldOfFractions(Rz);
assert Kz!F eq Evaluate(F0, 1/(Kz!z)-1)*(Kz!z)^6;
C := HyperellipticCurve(F);
J := Jacobian(C);
// Exact 11A certificate; A denotes [(1/2,1)-infinity].
A := J![C![1/2,1,1], C![1,0,0]];
M := 11*A;
u := z^2-(69/196)*z-9/392;
v := -(4085/2058)*z-233/1372;
assert M[1] eq u and M[2] eq v and M[3] eq 2;
assert (F-v^2) mod u eq 0;
trace_y := -Coefficient(u,1)*Coefficient(v,1)+2*Coefficient(v,0);
norm_y := Coefficient(v,0)^2
    -Coefficient(u,1)*Coefficient(v,0)*Coefficient(v,1)
    +Coefficient(u,0)*Coefficient(v,1)^2;
assert trace_y eq -139623/134456;
assert norm_y eq 7671/134456;
assert Valuation(trace_y,3) eq 1 and Valuation(norm_y,3) eq 1;
trace_y3 := trace_y^3-3*norm_y*trace_y;
assert Valuation(trace_y3,3) eq 3;
function Residue(q,n)
    R := Integers(n);
    return (R!Numerator(q))/(R!Denominator(q));
end function;
assert [Residue(Coefficient(u,i),81) : i in [2,1,0]] eq [1,48,63];
assert [Residue(3*Coefficient(v,i),81) : i in [1,0]] eq [14,75];
assert Residue(trace_y/3,3) eq 2;
assert Residue(norm_y/3,3) eq 2;
assert Residue(trace_y3/27,3) eq 1;
print "Exact 11A Mumford certificate:", M;
print "Trace and norm:", trace_y, norm_y;
print "Scaled residues S/3, N/3, Tr(y^3)/27 modulo 3:",
    Residue(trace_y/3,3), Residue(norm_y/3,3), Residue(trace_y3/27,3);
// First terms of the integral formal inverse and local primitives.
RY<Y> := PolynomialRing(Q);
inverse_z := Y^2-16*Y^4+560*Y^6;
assert (Evaluate(F,inverse_z)-Y^2) mod Y^8 eq 0;
omega0 := RY!(Derivative(inverse_z)/(2*Y));
omega1 := (inverse_z*omega0) mod Y^6;
assert Derivative(Y-(32/3)*Y^3+336*Y^5) eq omega0;
assert Derivative((1/3)*Y^3-(48/5)*Y^5) eq omega1;
assert -Residue(trace_y3/27,3)/Residue(trace_y/3,3) eq 1;
// Symbolic polynomial/root-map identities over Q(s).
Ks<s> := FunctionField(Q);
Rx<x> := PolynomialRing(Ks);
P := x^5+(2+s)*(3-2*s)*x^4
    +(1+10*s+4*s^2-14*s^3+s^4)*x^3
    -s*(1-14*s+4*s^2+10*s^3+s^4)*x^2
    -s^3*(1+2*s)*(-2+3*s)*x-s^5;
Kx := FieldOfFractions(Rx);
X := Kx!x;
den := (X+1)*(X+s^2);
Z := X*(s-1)^2/den;
W0 := (s-1)/(s+1);
Pstar := Rx!(-X^5*Evaluate(P,s^2/X)/s^5);
assert (Evaluate(F,Z)-W0^2)*den^5*(s+1)^2
    eq -(s-1)^2*(Kx!P)*(Kx!Pstar);
assert Evaluate(P,-1) eq 4*(s-1)^4*(s+1);
assert Evaluate(P,-s^2) eq -4*s^5*(s-1)^4*(s+1);
assert Evaluate(P,0) eq -s^5;
assert P-Pstar eq 8*(1-s^2)*x*(x+s)^3;
assert Evaluate(P,-s) eq -2*s^3*(s-1)^4;
// Distinct roots have the same Z-image only if their product is s^2.
// The reciprocal identities above exclude that case.
Rpair<xx,yy> := PolynomialRing(Ks,2);
assert xx*(yy+1)*(yy+s^2)-yy*(xx+1)*(xx+s^2)
    eq (xx-yy)*(s^2-xx*yy);
// The equal-y divisor curve and its elliptic quotient.
B<e,f,w> := PolynomialRing(Q,3);
H := 16*f^2-48*(e-1)^2*f+16*e*(e-1)^3-1;
G := w^2-16*f*((e-1)^3-(2*e-3)*f);
KB := FieldOfFractions(B);
Rzz<zz> := PolynomialRing(KB);
remainder := Evaluate(F,zz) mod (zz^2-e*zz+f);
assert Coefficient(remainder,1) eq -KB!H;
assert Coefficient(remainder,0)-(KB!w)^2 eq -KB!G;
V := 2*f-3*(e-1)^2;
assert V^2-((e-1)^3*(5*e-9)+1/4) eq H/4;
// Exact known rational lifts, and all six affine F_3 residue classes.
lifts := [<Q!3/2,Q!1/2,Q!1>, <Q!3/2,Q!1/2,Q!-1>,
          <Q!1,Q!1/4,Q!1>, <Q!1,Q!1/4,Q!-1>,
          <Q!1,Q!-1/4,Q!1>, <Q!1,Q!-1/4,Q!-1>];
for p in lifts do
    assert Evaluate(H,[p[1],p[2],p[3]]) eq 0;
    assert Evaluate(G,[p[1],p[2],p[3]]) eq 0;
end for;
assert F-1 eq -(z-1)*(2*z-1)^2*(4*z^2-4*z-1);
F3 := GF(3);
B3<e3,f3,w3> := PolynomialRing(F3,3);
H3 := 16*f3^2-48*(e3-1)^2*f3+16*e3*(e3-1)^3-1;
G3 := w3^2-16*f3*((e3-1)^3-(2*e3-3)*f3);
points3 := {<ee,ff,ww> : ee,ff,ww in F3 |
    Evaluate(H3,[ee,ff,ww]) eq 0 and Evaluate(G3,[ee,ff,ww]) eq 0};
expected3 := {<F3!0,F3!2,F3!1>, <F3!0,F3!2,F3!2>,
              <F3!1,F3!1,F3!1>, <F3!1,F3!1,F3!2>,
              <F3!1,F3!2,F3!1>, <F3!1,F3!2,F3!2>};
assert points3 eq expected3;
assert {<F3!p[1],F3!p[2],F3!p[3]> : p in lifts} eq points3;
He := Derivative(H3,1);
Hf := Derivative(H3,2);
for p in points3 do
    he := Evaluate(He,[p[1],p[2],p[3]]);
    hf := Evaluate(Hf,[p[1],p[2],p[3]]);
    assert hf ne 0 and p[3] ne 0;
    assert p[1]+he/hf ne 0;
end for;
print "Six residue classes on the equal-y divisor curve:", points3;
print "Six rational lifts:", lifts;
print "Nonvanishing differential checks: passed in all six disks.";
print "ALL STANDALONE CERTIFICATE ASSERTIONS PASSED.";
\end{code}

\noindent\textbf{Transcript.}
\begin{term}
Unconditional rank upper bound: 1
Complete genus-two rational-point list: { (1 : -8 : 1), (0 : -1 : 1), (1 : 8 :
1), (-1 : 0 : 1), (0 : 1 : 1), (1 : 0 : 0) }
Exact 11A Mumford certificate: (x^2 - 69/196*x - 9/392, -4085/2058*x - 233/1372,
2)
Trace and norm: -139623/134456 7671/134456
Scaled residues S/3, N/3, Tr(y^3)/27 modulo 3: 2 2 1
Six residue classes on the equal-y divisor curve: { <0, 2, 1>, <1, 2, 2>, <0, 2,
2>, <1, 1, 1>, <1, 2, 1>, <1, 1, 2> }
Six rational lifts: [ <3/2, 1/2, 1>, <3/2, 1/2, -1>, <1, 1/4, 1>, <1, 1/4, -1>,
<1, -1/4, 1>, <1, -1/4, -1> ]
Nonvanishing differential checks: passed in all six disks.
ALL STANDALONE CERTIFICATE ASSERTIONS PASSED.
\end{term}

The auxiliary primitive calculations in the script use the same normalization
$\omega_0=dz/(2y)$ and $\omega_1=z\,dz/(2y)$ as the main text.
The power-series assertions check their first coefficients; the all-orders precision estimates are proved in Lemma~\ref{lem:annihilator} and Proposition~\ref{prop:Cpoints}.

% ===================== AI assistance =====================

\section*{AI Assistance}
OpenAI GPT-6.1 Sol assisted with the development of the argument, symbolic calculations, verification code, and the preparation of the exposition. The author is responsible for the mathematical content and the final presentation.

% ===================== References =====================

\end{document}